\documentclass{IEEEtran}
\usepackage{hyperref}
\usepackage{times}
\usepackage{amsmath}
\usepackage{siunitx}

\usepackage{enumitem}
\usepackage{xcolor}
\usepackage[english]{babel}
\usepackage{amsfonts,amssymb,mathtools}
\usepackage{amstext}
\usepackage{upgreek}
\usepackage[tight,nice]{units}
\usepackage{algorithm2e}
\usepackage{eso-pic}

\IEEEpubid{}

\AddToShipoutPicture*{\footnotesize\sffamily\raisebox{1cm}{\hspace{1.65cm}\fbox{\parbox{\textwidth}{%
\copyright~2016 IEEE. Personal use of this material is permitted. Permission from IEEE must be obtained for all other uses, in any current or future media, including reprinting/republishing this material for advertising or promotional purposes, creating new collective works, for resale or redistribution to servers or lists, or reuse of any copyrighted component of this work in other works. DOI: \href{http://dx.doi.org/10.1109/URSI-EMTS.2016.7571330}{10.1109/URSI-EMTS.2016.7571330}.}}}}

\catcode`@=11
\newcommand{\frownfill}{$\scriptscriptstyle\m@th\mathord\frown\mkern-0.2mu%
     \cleaders\hbox{$\mkern-2mu\smash-\mkern-2mu$}\hfill
     \mkern-0.2mu$}
\newcommand{\bow}[1]{\vbox{\m@th\ialign{##\crcr%
      \frownfill\crcr\noalign{\kern-0.2\p@\nointerlineskip}%
      $\hfil\displaystyle{#1}\hfil$\crcr}}}
\newcommand{\bbow}[1]{\vbox{\m@th\ialign{##\crcr%
     \frownfill\crcr\noalign{\kern-0.7\p@\nointerlineskip}%
     \frownfill\crcr\noalign{\kern-0.3\p@\nointerlineskip}%
      $\hfil\displaystyle{#1}\hfil$\crcr}}}
\newcommand{\volt}[1]{\protect\bow{#1}}
\newcommand{\flux}[1]{\protect\bbow{#1}}
\newcommand{\fitmatrix}[1]{{\bf #1}}
\newcommand{\ve}{\volt{\rm\bf e}}
\newcommand{\vh}{\volt{\rm\bf h}}
\newcommand{\fj}{\flux{\rm\bf j}}
\newcommand{\MM}{\fitmatrix{M}}
\newcommand{\KK}{\fitmatrix{K}}
\newcommand{\A}{\fitmatrix{A}}
\newcommand{\Meps}{\MM_{\varepsilon}}
\newcommand{\Mmu}{\MM_{\mu}}

\newcommand{\CM}{\fitmatrix{C}}
\newcommand{\Cs}{\widetilde{\CM}}
\newcommand{\bn}{\mbox{\boldmath$n$}}
\newcommand{\be}{\mbox{\boldmath$e$}}
\newcommand{\bh}{\mbox{\boldmath$h$}}
\newcommand{\bd}{\mbox{\boldmath$d$}}
\newcommand{\bb}{\mbox{\boldmath$b$}}
\newcommand{\bj}{\mbox{\boldmath$j$}}
\newcommand{\CurlSymb}{\text{curl}}
\newcommand{\DivSymb}{\text{div}}
\newcommand{\Curl}[2][]  {\text{\CurlSymb}{#1}\,{#2}}
\newcommand{\Div}[2][]   {\text{\DivSymb}{#1}\,{#2}}

\begin{document}

\title{An Application of ParaExp to Electromagnetic Wave Problems} 

\author{\IEEEauthorblockN{    Melina Merkel\IEEEauthorrefmark{1}\IEEEauthorrefmark{2},
    Innocent Niyonzima\IEEEauthorrefmark{1}\IEEEauthorrefmark{2} and
    Sebastian Sch\"{o}ps\IEEEauthorrefmark{1}\IEEEauthorrefmark{2}}
  \medskip
  \IEEEauthorblockA{\IEEEauthorrefmark{1}Graduate School of Computational 
  Engineering (GSC CE), Technische Universit\"{a}t Darmstadt, Germany}
  \IEEEauthorblockA{\IEEEauthorrefmark{2}Institut f\"{u}r Theorie Elektromagnetischer Felder (TEMF), 
  Technische Universit\"{a}t Darmstadt, Germany}  
    e-mail: anna\_melina.merkel@stud.tu-darmstadt.de
  }

\maketitle

\begin{abstract}
    Recently, ParaExp was proposed for the time integration of hyperbolic problems. It splits the time interval of interest into
    sub-intervals and computes the solution on each sub-interval in parallel. The overall solution is decomposed into a particular
    solution defined on each sub-interval with zero initial conditions and a homogeneous solution propagated by the matrix exponential
    applied to the initial conditions. The efficiency of the method results from fast approximations of this matrix exponential using tools
    from linear algebra. This paper deals with the application of ParaExp to electromagnetic wave problems in time-domain. Numerical
    tests are carried out for an electric circuit and an electromagnetic wave problem discretized by the Finite Integration Technique.
\end{abstract}

\section{Introduction}

The simulation of high-frequency electromagnetic problems is often 
carried out in frequency domain. This choice is motivated by the linearity 
of the underlying governing equations and an implicit assumption that the source
signal can be decomposed in the Fourier basis.

However, the solution of problems in frequency domain may require the resolution of very
large linear systems of equations and this becomes particularly inconvenient for
broadband simulations involving multi-frequency sources such as Gaussian or pulsed
signals. The coupling with nonlinear time-dependent systems and the computation of
transients are other cases where time-domain simulations outperform frequency-domain
simulations.
On the other hand, the numerical complexity resulting from time-domain simulations may also become prohibitively expensive; time-domain parallelization is a promising solution alternative to domain decomposition in space.

The development and application of parallel-in-time methods dates back to more 
than 50 years \cite{nievergelt1964parallel}. These methods can be direct 
\cite{christlieb2010parallel,gander2013paraexp} or iterative 
\cite{lions2001parareal,minion2011hybrid}. 
They can also be well suited for small scale parallelization 
\cite{miranker1967parallel,womble1990time} or large 
parallelization \cite{gander2013paraexp,minion2011hybrid}.
Recently, the \emph{Parareal method} gained interest
\cite{lions2001parareal,gander2007analysis}. In its initial version, Parareal was 
developed for large scale parallelization of parabolic partial 
differential equations (PDEs). It involves the splitting of the time 
interval and the resolution of the governing ordinary differential equation (ODE) in parallel on each sub-intervals 
using a fine propagator which can be any classical time-stepper with a fine time grid. A coarse propagator distributes the initial conditions for each 
sub-interval during the Parareal iterations. It is typically obtained by a time stepper with a coarse grid on the entire time interval. Parareal iterates the 
resolution of both the coarse and the fine problems until convergence.

Most parallel-in-time methods fail for 
hyperbolic problems. In the case of Parareal, analysis has shown that it may lead to the \emph{beating} phenomenon depending on 
the structure of the system matrix obtained after space 
discretization \cite{farhat2006beat}. It may even become unstable if the eigenvalues 
of the matrix are pure imaginary which is the case in the 
presence of undamped electromagnetic waves.

\begin{figure}[t!]
  \centering
  \includegraphics[width=0.4\textwidth]{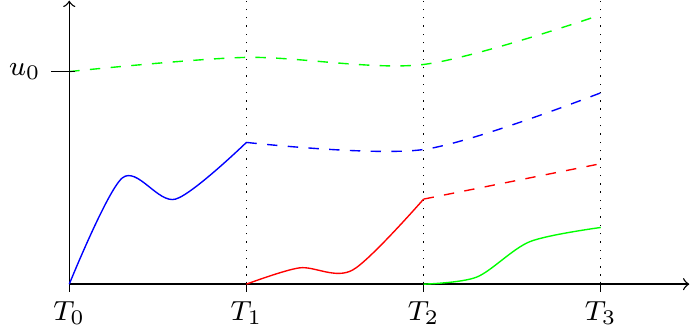}
  \vspace{-0.2cm}
  \caption{Schematic view of the decomposition of time and solution. Vertical dotted lines denote the sub-intervals, solid lines represent the solution of the inhomogeneous sub-problems and dashed lines represent the solution of the homogeneous sub-problems. Colors indicate the employed processors, cf. \cite{gander2013paraexp}}
  \label{fig:timedecomposition}
  \vspace{-0.2cm}
\end{figure}

In this paper we apply the \emph{ParaExp method} \cite{gander2013paraexp} for the parallelization 
of time-domain resolutions of hyperbolic equations that govern the electromagnetic 
wave problems. The method splits the time interval into sub-intervals and solves 
smaller problems on each sub-interval as visualized in Figure \ref{fig:timedecomposition}. Using the theory of linear 
ordinary differential equations, the total solution for each sub-interval is 
decomposed into particular solution with zero initial conditions and 
homogeneous solutions with initial conditions from previous intervals.
\IEEEpubidadjcol

The paper is organized as follows: in Section \ref{sec:maxwell-equations} we
introduce Maxwell's equations and derive the governing system of ODEs for the wave 
equation obtained by the Finite Integration Technique (FIT) or Finite Element Method (FEM). 
This system is then used in Section 
\ref{sec:paraexp} for the presentation of the ParaExp method following the lines of \cite{gander2013paraexp}. The mathematical 
framework is briefly sketched and the
details of the algorithm are discussed. Section \ref{sec:applications} deals with numerical 
examples. We consider two applications: a simple RLC circuit and an electromagnetic 
wave problem in an open wave guide. The accuracy of the method and the 
accuracy of the method in terms of the electromagnetic energy are investigated.

\section{Space and Time Discretization of Maxwell's equations}
\label{sec:maxwell-equations}

In an open, bounded domain $\Omega \subset \mathbb{R}^3$ 
and $t \in \mathcal{I}=(t_0,t_\mathrm{end}] \subset \mathbb{R}$,
the evolution of electromagnetic fields is governed by Maxwell's
Equations on $\Omega \times \mathcal{I}$, see e.g. \cite{jackson1999classical}:
\begin{align}
  \Curl[]{\bh}  = \bj +\partial_t \bd, \, \, \Curl[]{\be}  = -\partial_t \bb , 
  \, \,  \Div[]{\bd} = \rho, \, \, \Div[]{\bb}  = 0.
  \label{eq:maxwell-equations}
\end{align}
In presence of linear materials, these equations are completed by 
constitutive laws \cite{jackson1999classical}:
\begin{align}
    \bj=\sigma \be + \bj_s, \quad 
    \bd=\varepsilon\be, \quad 
    \bb=\mu\bh.
  \label{eq:const-laws-1}
\end{align}
In these equations, $\bh$ is the magnetic field [A/m], $\bb$ the magnetic flux
density [T], $\be$ the electric field [V/m], $\bd$ the electric flux density
[C/m$^2$], $\bj$, $\bj_{\rm d}=\partial_t \bd$ and $\bj_s$ are the Ohmic, displacement and electric source current densities [A/m$^2$], $\rho$ is the electric charge density [C/m$^3$]. 
The material properties $\sigma, \varepsilon$ and $\mu$ are the electric conductivity, 
the electric permittivity and the magnetic permeability, respectively.
In this paper, we consider electromagnetic wave propagation in
non-conducting media and which are free of charge, i.e., $\sigma \equiv 0$ 
and $\rho \equiv 0$. However, that the algorithms can be easily applied to 
the general case.

The space discretization of Maxwell's equations
\eqref{eq:maxwell-equations}-\eqref{eq:const-laws-1} 
using the Finite Integration Technique (FIT) \cite{Weiland_1977aa, Weiland_1996aa}
leads to the following initial
value problem (IVP)
\begin{align}
\MM \mathrm{d_t} \mathbf{u} + \KK \mathbf{u} = \bar{\mathbf{g}}(t) \quad \forall t \in \mathcal{I}, \quad \mathbf{u}(t_0)= \mathbf{u}_0.
\label{eq:FIT_ode_1}
\end{align}
with $\mathbf{u}^{\top} = [\vh^{\top},\ve^{\top}]$, $\bar{\mathbf{g}}(t)$ which is an excitation
and the matrices $\MM$ and $\KK$ given by:
\begin{align}
    \MM := 
    \begin{bmatrix}
        \Mmu^{\mathrm{FIT}} &  \boldsymbol{0} \\
        \boldsymbol{0}               & \Meps^{\mathrm{FIT}}  
    \end{bmatrix}
    \, , \, 
    \KK := 
    \begin{bmatrix}
        \boldsymbol{0}  & \CM \\
        -\Cs       & \boldsymbol{0}  
    \end{bmatrix}. \label{eq:fit_matrices}
\end{align}
$\Meps^{\mathrm{FIT}}$ and $\Mmu^{\mathrm{FIT}}$ are diagonal material matrices and $\CM$ and $\Cs$ are the discrete $\CurlSymb$ operators defined  
on the primal and dual grids, respectively; $\MM$ is invertible, thus \eqref{eq:FIT_ode_1} can be written as:
\begin{align}
\mathrm{d_t} \mathbf{u} = \A \mathbf{u}+{\mathbf{g}}(t) \quad \forall t \in \mathcal{I}, \quad \mathbf{u}(t_0)= \mathbf{u}_0
\label{eq:FIT_ode_2}
\end{align}
with $\A := -\MM^{-1} \KK$ and ${\mathbf{g}}(t) := \MM^{-1} \bar{\mathbf{g}}(t)$. Typically this system is time-integrated using the Leapfrog scheme.

Similarly, the use of the FE method applied to the $\be$-formulation leads 
to the following weak form: find $\be$ in an appropriate function space 
\cite{bossavit1998cem} such that
\begin{multline}
\left( \varepsilon \partial_{tt} \be, \be^{'}\right)_{\Omega}
+ \left( \mu^{-1} \Curl[]{\be}, \Curl[]{\be^{'}} \right)_{\Omega} 
\\
= -\left( \partial_{t} \bj_s, \be^{'}\right)_{\Omega_s}
+ \Big< (\mu^{-1}  \Curl[]{\be} \times \bn), \be^{'} \Big>_{\Gamma}
\label{eq:FEM_pde_1}
\end{multline}
holds for all test functions $\be^{'}$ in a space of test functions. 
In \eqref{eq:FEM_pde_1}, $\Gamma = \partial \Omega$ is the boundary of $\Omega$. 
Space discretization of \eqref{eq:FEM_pde_1} using the Galerkin approach 
leads to the system of ODEs:
\begin{equation}
\begin{aligned}
&\Meps^{\mathrm{FE}} \mathrm{d_{tt}} \be + 
\KK_{\mu^{-1}}^{\mathrm{FE}} \be 
= 
\bar{\mathbf{g}}^{\mathrm{FE}}(t) \quad \forall t \in \mathcal{I},
\\
&\be(t_0)= \be_0,\quad\mathrm{d}_t \be(t_0)= \be_1.
\label{eq:FEM_ode_2}
\end{aligned}
\end{equation}
where the matrices $\Meps^{\mathrm{FE}}$ and $\KK_{\mu^{-1}}^{\mathrm{FE}}$ in \eqref{eq:FEM_ode_2} are FEM material 
matrices obtained from the discretization of the bilinear forms in 
\eqref{eq:FEM_pde_1} and $\bar{\mathbf{g}}^{\mathrm{FE}}(t)$ is the source term.
Equation \eqref{eq:FEM_ode_2} can be recast in a first order system similar to \eqref{eq:FIT_ode_1}.

\section{The ParaExp Algorithm}
\label{sec:paraexp}

In this section we develop based on \cite{gander2013paraexp} the ideas of the ParaExp method for the system of ODEs
in the form \eqref{eq:FIT_ode_2}
\begin{equation}
\begin{aligned}
\mathrm{d_t} \mathbf{u} =  \mathbf{\A} \mathbf{u}+\mathbf{g}(t) \quad \forall t \in \mathcal{I}, \quad \mathbf{u}(t_0)= \mathbf{u}_0.
\label{eq:FIT_ode_3}
\end{aligned}
\end{equation}

Applying the method of variation of constants to equation \eqref{eq:FIT_ode_3}
leads to the solution
\begin{equation}
\mathbf{u}(t) = \exp{(t \A}) \mathbf{u}_0 + \int_0^t \exp{((t-\tau) 
\A ) }\mathbf{g}(\tau) \mathrm{d} \tau
\label{eq:solution-variation-of-constants}
\end{equation}
where $\exp{(t \A}) \mathbf{u}_0$ is the homogeneous solution due to 
initial conditions and the convolution product is the particular solution resulting 
from the presence of the source term $\mathbf{g}(t)$. The last term of 
\eqref{eq:solution-variation-of-constants} is more difficult to compute than the first
one. However, thanks to the linearity of equation \eqref{eq:FIT_ode_3} and the superposition 
principle, $\mathbf{u}(t)$ can be written as $\mathbf{u}(t) = \mathbf{v}(t) + \mathbf{w}(t)$ where the particular
solution $\mathbf{v}(t)$ is governed by
\begin{equation}
\mathrm{d_t} \mathbf{v} = \A \mathbf{v}+\mathbf{g}(t) \quad \forall t \in \mathcal{I}, \quad \mathbf{v}(t_0)= \boldsymbol{0}
\label{eq:solution-particular}
\end{equation}
and the homogeneous solution $\mathbf{w}(t)$ is governed by
\begin{equation}
\mathrm{d_t} \mathbf{w} =\A \mathbf{w}+ \boldsymbol{0} \quad \forall t \in \mathcal{I}, \quad \mathbf{w}(t_0)= \mathbf{u}_0.
\label{eq:solution-homogeneous}
\end{equation}
The ParaExp method takes advantage of this decomposition.
The time interval (0, T] is partitioned into sub-intervals 
$\mathcal{I}_j = (T_{\mathrm{j-1}}, T_{\mathrm{j}}]$
with $j = 1, 2, ..., p$,  $t_0 = T_0 < T_1 < T_2 < ... < T_p = t_{\mathrm{end}}$
and $p$ the number of CPUs. The following solutions are then computed on each CPU:
\begin{enumerate}[label=(\bfseries\alph*)]
\item a particular solution $\mathbf{v}_j(t)$ governed by: \label{enum:paraexpA}
\begin{equation}
\label{eq:solution-particular-paraexp}
\begin{aligned}
\mathrm{d_t} \mathbf{v}_j= \A \mathbf{v}_j +\mathbf{g}(t) \quad 
& \forall t \in \mathcal{I}_j,
\\
& \mathbf{v}_j(T_{j-1})= \boldsymbol{0},
\end{aligned}
\end{equation}
\item a homogeneous solution $\mathbf{w}_j(t)$ governed by: \label{enum:paraexpB}
\begin{equation}
\label{eq:solution-homogeneous-paraexp}
\begin{aligned}
\mathrm{d_t} \mathbf{w}_j = \A \mathbf{w}_j +\boldsymbol{0} \quad
&\forall t \in (T_{j-1}, T], 
\\
&\mathbf{w}_j(T_{j-1})= \mathbf{v}_{\mathrm{j-1}}(T_{j-1}).
\end{aligned}
\end{equation}
\end{enumerate}
Problems \eqref{eq:solution-particular-paraexp} can be solved simultaneously 
in parallel using a time stepping method as no initial conditions need to be provided. 
Problems \eqref{eq:solution-homogeneous-paraexp} are solved on (possibly bigger) 
time intervals $(T_{j-1}, T]$ and they yield exponential solutions 
$\mathbf{w}_j(t) = \exp{(t \A)} \mathbf{v}_{\mathrm{j-1}}(T_{j-1})$
where the initial condition is the final solution $\mathbf{v}_{\mathrm{j-1}}(T_{j-1})$. 
It is therefore highly recommended to compute $\mathbf{v}_{\mathrm{j-1}}(t)$ and 
$\mathbf{w}_{\mathrm{j}}(t)$ on the same CPU to avoid communicational costs.
Using the superposition principle, the total solution can be expanded as:
\begin{equation}
\mathbf{u}(t) = \mathbf{v}_j(t) + \sum_{i = 1}^{j}\mathbf{w}_i(t) \, \, \textrm{ with } j \textrm{ such that } 
t\in \mathcal{I}_j.
\label{eq:solution-particular-para}
\end{equation}
Figure \ref{fig:timedecomposition} shows the time decomposition of IVP into 
particular solutions (solid lines) and homogeneous solutions (dashed lines)
for a case with 3 CPUs.
Steps \ref{enum:paraexpA} and \ref{enum:paraexpB} of the ParaExp method are described in Algorithm \ref{alg:ParaExp}.

\begin{center}
\vspace*{-0.5em}
\begin{algorithm}
\SetKwInput{Input}{Input}
\SetKwInput{Output}{Output}
\LinesNumbered
\DontPrintSemicolon
\Input{system matrix $\mathbf{A}$, source term $\mathbf{g}(t)$, initial value $\mathbf{u}_0$, time interval $\mathcal{I}$, number of processors $p$}
\Output{solution $\mathbf{u}(t)$}
\Begin
{   set $T_j$ (partition $\mathcal{I}$ into intervals $\mathcal{I}_j$, $j=0,...,p$),\\
    \emph{\# begin the parallel loop (index $j$)}\\
    \For{$(j \gets 1$ \KwTo $p)$}{
	$\mathbf{v}_j \gets$ solve $\mathbf{v}_j'(t)=\mathbf{A}\mathbf{v}_j(t)+\mathbf{g}(t)$,~$\mathbf{v}_j(T_{j-1})=0$, $t \in \mathcal{I}_j$ using a time stepper\\
	\eIf{$j\neq p$}{
        $\mathbf{w}_{j+1} \gets \exp(\mathbf{A}(t-T_{j}))\mathbf{v}_j(T_j)$ for all $t \in (T_j,T_p]$
        }{
        $\mathbf{w}_{1} \gets \exp(\mathbf{A}(t-T_{0}))\mathbf{u}_0$ for all $t \in (T_0,T_p]$
        }
    }
    \For{$(j \gets 1$ \KwTo $p)$}{
    $\mathbf{u}(t) \gets \mathbf{v}_j(t)+\sum\limits_{i=1}^j\mathbf{w}_i(t)$, for all $t\in \mathcal{I}_j$
    }
}
\vspace{1em}
\caption{Pseudocode for the ParaExp Algorithm}
\label{alg:ParaExp}
\end{algorithm}
\vspace*{-1em}
\end{center} 
{A critical point of the method is the computation of the matrix exponential. 
A straight forward evaluation of the definition is not feasible. Instead, efficient approximations
of the action of the matrix exponential should be used, e.g., based on Krylov subspaces 
\cite{gander2013paraexp} or the approximation of the action of the matrix exponential 
as proposed in \cite{alhi11}, where for given integers $s$ and $m$ and for an arbitrary vector $\mathbf{b}$ the relation
\begin{align*}
\exp(t\A)\mathbf{b} = \left(\exp\!\left(s^{-1}t\A\right)\right)^s\mathbf{b}
\end{align*}
is exploited to build a recurrence 
\begin{align*}
\mathbf{u}_{i+1}=\mathbf{r}_m\left(s^{-1}t\A\right)\mathbf{b}_i, \;\;\; i=0,\dots ,s-1, \;\;\; \mathbf{b}_0=\mathbf{b},
\end{align*}
with the truncated Taylor series of order $m$ of the matrix exponential
\begin{align*}
 \mathbf{r}_m\left(s^{-1}t\A\right) = \sum\limits_{j=0}^m\frac{\left(s^{-1}t\A\right)^j}{j!}.
\end{align*}} 

\section{Preliminary results}
\label{sec:applications}
\begin{figure}[h!]
  \centering
  \includegraphics[width=0.4\textwidth]{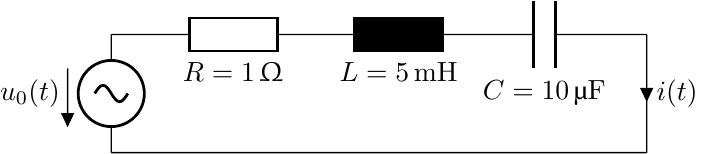}
  \caption{RLC test circuit}
 \label{fig:rlc}
\end{figure}

\begin{figure}[h!]
  \centering
  \includegraphics[width=0.49\textwidth]{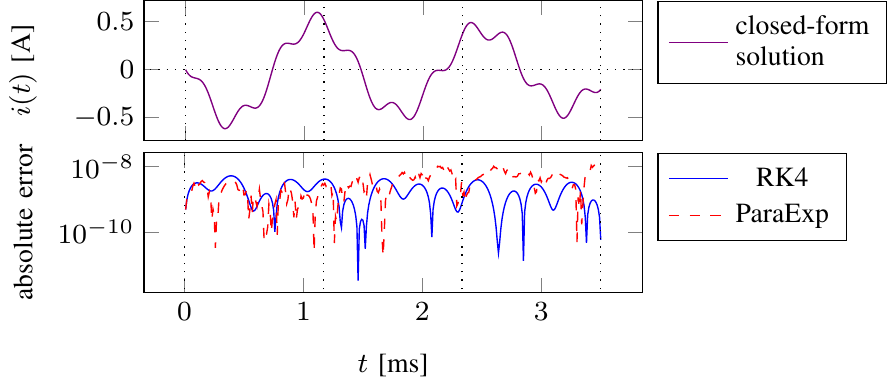}
  \vspace{-0.7cm}
  \caption{Closed form solution for $i(t)$ and absolute error of RK4 and ParaExp}
  \label{fig:rlcCurrent}
\end{figure}

\begin{figure}[h!]
  \centering
  \includegraphics[width=0.47\textwidth]{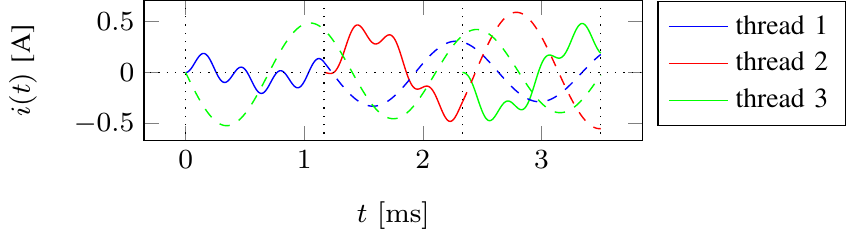}
  \vspace{-0.3cm}
  \caption{Time decomposition of the problem into three inhomogeneous problems with zero initial value (solid lines) and three homogeneous problems (dashed lines)}
 \label{fig:rlcDecomp}
\end{figure}

As a proof-of-concept, the ParaExp algorithm is implemented in Octave \cite{octave}, using the explicit Runge-Kutta method of order 4 (RK4) provided by OdePkg \cite{odepkg} as a time stepper, while no approximation of the matrix exponential was used. The RLC circuit (Fig. \ref{fig:rlc}) with known closed-form solution is constructed as a first test case.
The signal of the voltage source is $u_0(t)=U_0\sin\left(\omega_0 t\right)$ with $U_0=\SI{10}{\volt}$ and $\omega_0=2000\pi^2\si[per-mode=fraction]{\per\second}$. The differential equation of this problem is given by
\begin{align}
  L\mathrm{d_{tt}^2}i(t)+R \mathrm{d_{t}}i(t)+C^{-1} i(t)=U_0\omega_0\cos\left(\omega_0 t\right) \label{eq:RLCdiffEq}\\
  i(0)=0,\qquad  \mathrm{d_{t}}i(0)=-{U_{L0}}{L}^{-1}
\end{align}
with $U_{L0}=\SI{12}{\volt}$.
Both, RK4 on the whole time interval and ParaExp with three parallel threads, are applied to equation \eqref{eq:RLCdiffEq}. ParaExp executes the RK4 time stepper in each thread and the constant time step size has been chosen in all cases as $\Delta t=\SI{1e-5}{s}$. The results are compared to the closed-form solution in Fig. \ref{fig:rlcCurrent}.
As expected, the error of the ParaExp algorithm is of the same order of magnitude as the one of the traditional RK4 time stepper. An illustration of the decomposition of the problem on three parallel threads is given in Figure~\ref{fig:rlcDecomp}.

\begin{figure}[h!]
  \centering
  \includegraphics[width=0.4\textwidth]{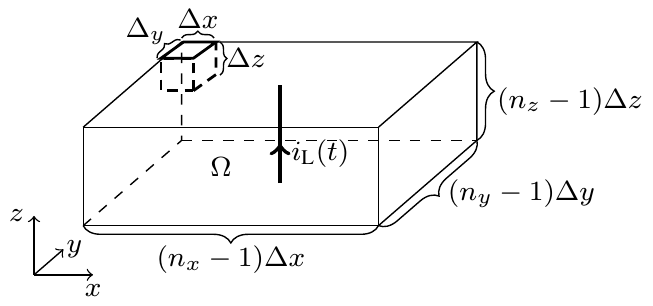}
  \caption{Domain $\Omega$ of the wave problem with a hexahedral mesh}
  \label{fig:waveProbDef}
\end{figure}

To test the algorithm on a more complex example a 2D cylindrical wave is simulated. The excitation is a line current in $z$-direction in the center of the domain $\Omega$ as shown in Figure~\ref{fig:waveProbDef}. 
The discretization is obtained by FIT. A PEC boundary is assumed on the whole boundary $\partial \Omega$. The parameters of the discretization are $\Delta x=\Delta y=\Delta z=\SI{1}{\meter}$ and $n_x=n_y=21,\: n_z=2$ and the domain is filled with vacuum.
The line current is a Gaussian function given by
\begin{align}
 i_\text{L}(t)=i_\text{max}e^{-4\left(\frac{t-\sigma_t}{\sigma_t}\right)^2} \label{eq:line_current}
\end{align}
with $i_\text{max}=\SI{1}{\A}$ and $\sigma_t=\SI{2e-8}{\s}$.\\
 
The differential equation of this problem is given by \eqref{eq:FIT_ode_1} and \eqref{eq:fit_matrices} with $\bar{\mathbf{g}}(t) = -[\boldsymbol{0},\fj]^{\top}$, with $\fj$ being the discretized line current \eqref{eq:line_current}.

\begin{figure}[h]
  \centering
  \includegraphics[width=0.47\textwidth]{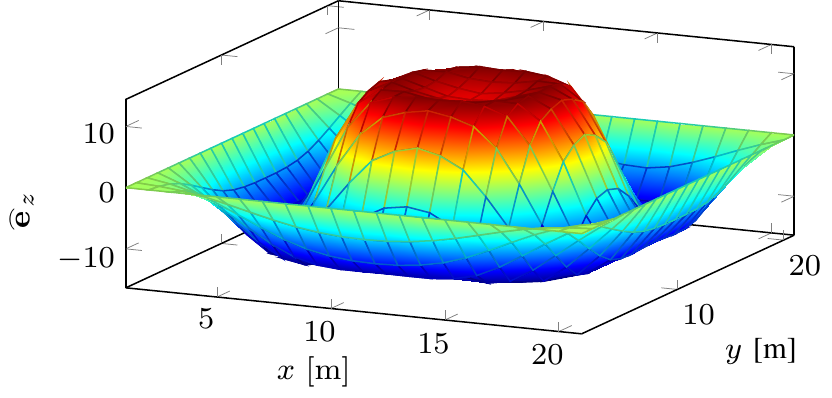}
  \caption{The $\ve_z$ component of the wave at $t=\SI{4.4e-8}{\s}$}
 \label{fig:wave}
\end{figure}

\begin{figure}[h!]
  \centering
  \includegraphics[width=0.47\textwidth]{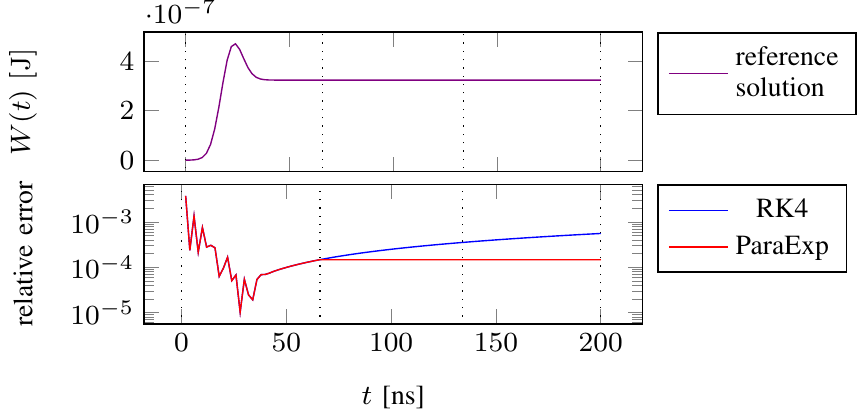}
  \vspace{-0.2cm}
  \caption{Reference solution (RK4 with a relative tolerance of \SI{1e-10}{}) for $W(t)$ and relative error of RK4 (with fixed time steps) and ParaExp}
 \label{fig:waveEnergy}
\end{figure}

 The $\ve_z$ component of the calculated wave can be seen in Figure \ref{fig:wave} at $t=\SI{4.4e-8}{\s}$.
RK4 on the whole time interval and ParaExp (without approximation of the matrix exponential) on three threads are used to solve \eqref{eq:FIT_ode_1}. Both approaches use a constant time step size $\Delta t=\SI{2e-9}{s}$. With the obtained results, the energy in the domain is calculated using
\begin{align}
  W=\frac{1}{2}\left(\ve^{\!\top}\Meps\ve+\vh^{\!\top}\Mmu\vh\right).  
\end{align}
The results are compared to a more accurate reference solution (i.e. the solution of the Runge-Kutta solver with a relative error tolerance of \SI{1e-10}{}).
 The comparison of the relative error of ParaExp and of the traditional RK4 method can be seen in Figure \ref{fig:waveEnergy}.
 On the first interval both approaches coincide, afterwards the exponential is more accurate than the time stepper. However, in the current implementation this comes with higher costs since no approximation of the matrix exponential is used. 
\section{Outlook}
In the extended paper, we will discuss the application of the ParaExp method in 
more detail. The approximation of the action of the matrix exponential, the utilization of the Leapfrog scheme
within ParaExp, and the numerical costs of the various methods will be investigated and compared. Finally the computational efficiency will be discussed using more complex examples.
~\\

\section*{Acknowledgment}
The authors would like to thank Timo Euler, CST AG for the fruitful discussions on time domain simulations.

This work was supported by the German Funding Agency (DFG) by the grant 
`Parallel and Explicit Methods for the Eddy Current Problem' (SCHO-1562/1-1), 
the 'Excellence Initiative' of the German Federal and State Governments and the 
Graduate School CE at Technische Universit\"at Darmstadt.
~\\


\begin{thebibliography}{10}
\providecommand{\url}[1]{#1}
\csname url@samestyle\endcsname
\providecommand{\newblock}{\relax}
\providecommand{\bibinfo}[2]{#2}
\providecommand{\BIBentrySTDinterwordspacing}{\spaceskip=0pt\relax}
\providecommand{\BIBentryALTinterwordstretchfactor}{4}
\providecommand{\BIBentryALTinterwordspacing}{\spaceskip=\fontdimen2\font plus
\BIBentryALTinterwordstretchfactor\fontdimen3\font minus
  \fontdimen4\font\relax}
\providecommand{\BIBforeignlanguage}[2]{{%
\expandafter\ifx\csname l@#1\endcsname\relax
\typeout{** WARNING: IEEEtran.bst: No hyphenation pattern has been}%
\typeout{** loaded for the language `#1'. Using the pattern for}%
\typeout{** the default language instead.}%
\else
\language=\csname l@#1\endcsname
\fi
#2}}
\providecommand{\BIBdecl}{\relax}
\BIBdecl

\bibitem{nievergelt1964parallel}
J.~Nievergelt, ``Parallel methods for integrating ordinary differential
  equations,'' \emph{Communications of the ACM}, vol.~7, no.~12, pp. 731--733,
  1964.

\bibitem{christlieb2010parallel}
A.~J. Christlieb, C.~B. Macdonald, and B.~W. Ong, ``Parallel high-order
  integrators,'' \emph{SIAM Journal on Scientific Computing}, vol.~32, no.~2,
  pp. 818--835, 2010.

\bibitem{gander2013paraexp}
M.~J.~Gander and S.~G\"uttel, ``{P}ara{E}xp: A parallel integrator for linear
  initial-value problems,'' \emph{SIAM Journal on Scientific Computing},
  vol.~35, no.~2, pp. C123--C142, 2013.

\bibitem{lions2001parareal}
J.-L. Lions, Y.~Maday, and G.~Turinici, ``A `{P}arareal' in time discretization
  of {PDE}s,'' \emph{Comptes Rendus de l'Academie des Sciences Series I
  Mathematics}, vol. 332, no.~7, pp. 661--668, 2001.

\bibitem{minion2011hybrid}
M.~Minion, ``A hybrid {P}arareal spectral deferred corrections method,''
  \emph{Communications in Applied Mathematics and Computational Science},
  vol.~5, no.~2, pp. 265--301, 2011.

\bibitem{miranker1967parallel}
W.~L. Miranker and W.~Liniger, ``Parallel methods for the numerical integration
  of ordinary differential equations,'' \emph{Mathematics of Computation},
  vol.~21, no.~99, pp. 303--320, 1967.

\bibitem{womble1990time}
D.~E. Womble, ``A time-stepping algorithm for parallel computers,'' \emph{SIAM
  Journal on Scientific and Statistical Computing}, vol.~11, no.~5, pp.
  824--837, 1990.

\bibitem{gander2007analysis}
M.~J. Gander and S.~Vandewalle, ``Analysis of the {P}arareal time-parallel
  time-integration method,'' \emph{SIAM Journal on Scientific Computing},
  vol.~29, no.~2, pp. 556--578, 2007.

\bibitem{farhat2006beat}
C.~Farhat, J.~Cortial, C.~Dastillung, and H.~Bavestrello, ``{T}ime-parallel
  implicit integrators for the near-real-time prediction of linear structural
  dynamic responses,'' \emph{International Journal for Numerical Methods in
  Engineering}, vol.~67, no.~5, pp. 697--724, 2006.

\bibitem{jackson1999classical}
J.~D. Jackson, \emph{Classical Electrodynamics}.\hskip 1em plus 0.5em minus
  0.4em\relax Wiley, 1999.

\bibitem{Weiland_1977aa}
T.~Weiland, ``A discretization model for the solution of {Maxwell}'s equations
  for six-component fields,'' \emph{International Journal of Electronics and
  Communications}, vol.~31, pp. 116--120, 1977.

\bibitem{Weiland_1996aa}
------, ``Time domain electromagnetic field computation with finite difference
  methods,'' \emph{International Journal of Numerical Modelling: Electronic
  Networks, Devices and Fields}, vol.~9, no.~4, pp. 295--319, 1996.

\bibitem{bossavit1998cem}
A.~Bossavit, \emph{Computational Electromagnetism. {V}ariational Formulations,
  Complementarity, Edge Elements}.\hskip 1em plus 0.5em minus 0.4em\relax
  Academic Press, 1998.

\bibitem{alhi11}
A.~H. Al-Mohy and N.~J. Higham, ``Computing the action of the matrix
  exponential, with an application to exponential integrators,'' \emph{SIAM
  Journal on Scientific Computing}, vol.~33, no.~2, pp. 488--511, 2011.

\bibitem{octave}
\BIBentryALTinterwordspacing
J.~W. Eaton \emph{et~al.}, ``{GNU} {O}ctave.'' [Online]. Available:
  \url{http://www.octave.org}
\BIBentrySTDinterwordspacing

\bibitem{odepkg}
\BIBentryALTinterwordspacing
T.~Treichl and J.~Corno, \emph{{OdePkg}, A Package for Solving Differential
  Equations with {O}ctave}.\hskip 1em plus 0.5em minus 0.4em\relax Free
  Software Foundation. [Online]. Available:
  \url{http://octave.sourceforge.net/odepkg/}
\BIBentrySTDinterwordspacing

\end{thebibliography}
\end{document}